\documentclass[11pt]{article}
\usepackage{dsfont, amssymb,amsmath,amscd,latexsym, amsthm, amsxtra,amsfonts}
\usepackage[all]{xy}
\newtheorem{Them}{Theorem}[section]
\newtheorem{Prop}{Proposition}[section]

\newtheorem{Lemma}{Lemma}[section]

\begin{document}
\makeatletter
\def\@setauthors{%
\begingroup
\def\thanks{\protect\thanks@warning}%
\trivlist \centering\footnotesize \@topsep30\p@\relax
\advance\@topsep by -\baselineskip
\item\relax
\author@andify\authors
\def\\{\protect\linebreak}%
{\authors}%
\ifx\@empty\contribs \else ,\penalty-3 \space \@setcontribs
\@closetoccontribs \fi
\endtrivlist
\endgroup } \makeatother
 \baselineskip 15pt
\title{{\Large {\bf Large deviations for multidimensional SDEs with
reflection}}
\thanks{ This work is
 supported by NSFC and  SRF for ROCS, SEM}}
\author{ Zongxia \ Liang \thanks{ email: zliang@math.tsinghua.edu.cn}
\\ \small Department of Mathematical Sciences,
\  Tsinghua University,\\
\small  Beijing 100084, People's Republic of China }
\date{}
 \maketitle
\begin{abstract}\noindent
    The  large deviations  principles are  established  for a class of multidimensional
degenerate stochastic differential equations  with  reflecting
boundary conditions. The results include two cases where the initial
conditions are  adapted and anticipated. \vskip 10pt \noindent {\bf
MSC}(2000): Primary 60F10, 60H10 ; Secondary  60J50, 60J60. \vskip
10pt \noindent
\end{abstract}
 \setcounter{equation}{0}
\section{ {\small {\bf Introduction and main results}}}

\noindent Let $\mathcal{O}$ be a  smooth  bounded  open set in $
\Re^d $. $ {\bf n}(x) $ denotes the cone of unit outward normal
vectors to $
\partial \mathcal{O} $ at $x$, that is,
\begin{eqnarray}
&(i) & \exists \ C_0 \geq 0,  \forall x\in \partial \mathcal{O},
\forall x' \in \bar{\mathcal{O}} ,  \exists \  k \in {\bf
n}(x)\nonumber\\
&& \Longrightarrow  (x-x', k) + C_0|x-x'|^2\geq 0
,\\
& (ii) &  \forall  x\in  \partial \mathcal{O}, \mbox{if} \  \exists
C \geq 0, \exists \ k \in \Re^d  , \forall  x' \in \bar{
\mathcal{O}}
, \nonumber\\
&& (x-x', k ) + C|x-x'|^2\geq 0 , \Longrightarrow  k =\theta {\bf
n}(x)\\
&& \mbox{ for some  $  \theta \geq 0$. Moreover, we assume that  }\nonumber\\
& (iii) &  \exists n \geq 1, \exists \alpha >0, \exists R>0; \exists
a_1, \cdots, a_n\in \Re^d, |a_i|=1, \forall i, \nonumber\\
&& \exists x_1, \cdots, x_n \in \partial \mathcal{O}: \partial
\mathcal{O}\subset \cup^n_{i=1} B(x_i, R),\nonumber\\
&& \forall i, \forall x \in \partial \mathcal{O}\cap B(x_i, 2R),
\forall \xi \in {\bf n}(x), \ \ (\xi, a_i)\geq \alpha >0,
\end{eqnarray}
 where $ \partial \mathcal{O} $ denotes
the boundary of $\mathcal{O}$, $\bar{\mathcal{O}} $ denotes the
closure of $\mathcal{O}$, $ B(x,r)$ denotes the ball of  radius $r$
at $x$. We assume that $ B_t $ is an  $ \Re^d $- valued
$\mathcal{F}_t $- Brownian motion on  a stochastic basis $(\Omega,
{\mathcal F},\{{\mathcal F}_t\}_{t\in [0,1]}, {\bf P})$ satisfying
the usual assumptions. For $\varepsilon >0  $ we consider  the
following perturbed stochastic differential equations on domain $
\mathcal{O} $ with reflecting boundary conditions,
\begin{eqnarray} \left \{\begin{array}{l}
X^\varepsilon_t(x)=x + \int^t_0 b(X^\varepsilon_s(x))ds
+\sqrt{\varepsilon}\int^t_0 \sigma (X^\varepsilon_s(x))dB_s
-L^\varepsilon_t(x), \\
L^\varepsilon_t(x)=\int^t_0\xi(X^\varepsilon_s(x))d
|L^\varepsilon(x)|_s,\\
|L^\varepsilon(x)|_t=\int^t_0I_{\{s:X^\varepsilon_s(x)\in \partial
\mathcal{O} \}}d|L^\varepsilon(x)|_s
\end{array} \right.
\end{eqnarray}
for $\forall \ t \in [0,1]$ and  $ x \in \bar{\mathcal{O}}$, where
$b: \Re^d\mapsto \Re^d $ and $ \sigma: \Re^d \mapsto \Re^d \times
\Re^d $ are continuous functions, ${\bf \xi }(X^\varepsilon_s(x))
\in {\bf n}(X^\varepsilon_s(x)) $,  the $ |L^\varepsilon(x)|_t $
denotes the total variation of $ L^\varepsilon_t(x)$ on $[0, t]$. A
pair  $ (X^\varepsilon_t(x), L^\varepsilon_t(x))$ of continuous
processes  is called a solution to equations (1.4) if there exists a
measurable set $\tilde{ \Omega } $ with ${\bf P}( \tilde{ \Omega
})=1$ such that for each $\omega \in \tilde{ \Omega } $ (i) for each
$ x \in \bar{\mathcal{O}}$ the function $s \mapsto$ $
L^\varepsilon_s(x)$ with  values in $\Re^d $  has bounded variation
on  interval $ [0, 1]$ and $ L_0^x=0 $; (ii) for all $t\geq 0 $,
$X^\varepsilon_t(x)\in \bar{\mathcal{O}}$ and $ (X^\varepsilon_t(x),
L^\varepsilon_t(x))$ satisfies Eq.(1.4).\vskip 19pt

We first recall that the existence and uniqueness results on strong
solutions to Eq.(1.4) were  studied by Skorohod\cite{s1},
Tanaka\cite{s2}, Lions and Sznitman\cite{s3},  Saisho\cite{s4}, and
other authors. When the initial value $x$ is replaced by an
arbitrary random variable $\{ X^\varepsilon_0, \varepsilon>0 \}$,
the author and Zhang (see \cite{s5, s6}) obtained recently existence
of strong solution to Eq.(1.4), and we proved that the composition $
(X^\varepsilon_t( X^\varepsilon_0 ),
L^\varepsilon_t(X^\varepsilon_0))$ of stochastic processes
 $ (X^\varepsilon_t(x), L^\varepsilon_t(x))$  and $ X^\varepsilon_0  $
 is just a solution of Eq.(1.4) corresponding to the initial
  $ X^\varepsilon_0$ under little regularity  conditions on $b$,
  $ \sigma $ and shape of domain $ \mathcal{O}$. \\

The goal of  this paper will be  two-fold. One is to establish the
large deviation principle  for  $\{X^\varepsilon_t(x): \varepsilon
>0 \}$. We noted that this problem has been solved by Anderson and
Orey \cite{s7}, C\'{e}pa\cite{s8} and other authors. However, this
result only deals with the case where the $\sigma \sigma^T $ is
uniformly definite, i.e., $\sigma $ is non-degenerate,  because
their proofs  heavily depend  on one dimensional stochastic
differential equations with reflection in which the reflection has
explicit representation formula. Here  we will remove these
restrictions  and give detailed proof of the problem in case where $
\sigma $ is degenerate by a approach used by the author\cite{s9},
which is somewhat different from that of \cite{s7,s8}. The other one
is to prove that the non-adapted solution $ \{X^\varepsilon_t(
X^\varepsilon_0 )\} $ also satisfies the large deviation principle
under  some hypotheses on the family  $\{ X^\varepsilon_0,\
\varepsilon>0\}$ via some results on adapted solution
$\{X^\varepsilon_t(x),\ \varepsilon >0 \}$.\vskip 19pt

 To state our result more precisely, we introduce the following
 skeleton  equation associated with (1.4),
\begin{eqnarray} \left \{\begin{array}{l}
z^\psi_t(x)=x + \int^t_0 b(z^\psi_s(x))ds +\int^t_0 \sigma
(z^\psi_s(x))\psi(s)ds-k^\psi_t(x), \\
k^\psi_t(x)=\int^t_0\xi(z^\psi_s(x))d|k^\psi(x)|_s,\\
|k^\psi(x)|_t=\int^t_0I_{\{s:z^\psi_s(x)\in \partial \mathcal{O}
\}}d|k^\psi(x)|_s
\end{array} \right.
\end{eqnarray}
where $\psi\in L^2([0,1]; \Re^d)$, the $ |k^\psi(x)|_t $ denotes the
total variation of bounded variation function $ k^\psi_t(x)$ on $[0,
t]$ with $ k^\psi_0(x)=0$ . A pair $ (z^\psi_t(x), k^\psi_t(x))$ of
functions is a solution of (1.5) means that $(z^\psi_t(x),
k^\psi_t(x))\in \bar{\mathcal{O}}\times \Re^d $
satisfies Eq.(1.5).\\

Let $E=C([0,1];\bar{\mathcal{O}})$. $\|\cdot \|_{E} $ denotes the
uniform norm on $E$. We define the rate functions $I_1$, $I^x_2$ and
$I^x$ by
\begin{eqnarray}  I_1(g)= \left  \{\begin{array}{l}
\frac{1}{2}\int^1_0|g(s)|^2ds,\  g\in L^2([0,1]; \Re^d  ),\\
+\infty,  \ \ \  \   \qquad  \quad \mbox{otherwise, }
\end{array}\quad\quad\  \right.
\end{eqnarray}
\begin{eqnarray}
 &&I^x_2(f)=\inf\{I_1(  \psi): f=z^\psi(x)\ \mbox{  is a solution of
 Eq.(1.5)} \},\\
 &&I^x(f)=\limsup_{y\rightarrow x}\{I^y_2(f)\}
\end{eqnarray}
for $x\in \bar{\mathcal{O}} $  and $ f\in E$. We set $ \inf
\varnothing = +\infty $ by convention. The main results of this
paper are  the following.
\begin{Them}
Assume that $\mathcal{O} $ is a smooth bounded open set in $\Re^d$
satisfying (1.1)-(1.3) and there exists a function $\phi\in
{\mathcal{C}}_b^2(\Re^d ) $ such that
\begin{eqnarray}
\exists \alpha >0,\ \forall x \in \partial \mathcal{O},\ \forall
\zeta \in {\bf n}(x),\ (\triangledown \phi(x), \zeta ) \leq -\alpha
C_0.
\end{eqnarray}
Let $b$ and $\sigma $ be bounded and uniformly Lipschitz. $\{
X^\varepsilon_t(x), \varepsilon >0 \}$  is the uniqueness solution
of Eq.(1.4). Then  $\{ X^\varepsilon(x), \varepsilon >0 \}$
satisfies  large deviation principle on $E$ with
 good rate function $ I^x_2$ given by (1.7). In other words, for any
 open set $G$ and any closed set $F$ of $E$, we have
\begin{eqnarray}
&& \liminf_{\varepsilon \rightarrow 0}\varepsilon \log {\bf P}\{
X^\varepsilon(x)\in G\}\geq -\inf_{f\in G}\{I^x_2(f)\},\\
&& \limsup_{\varepsilon \rightarrow 0}\varepsilon \log {\bf P}\{
X^\varepsilon(x)\in F\}\leq  -\inf_{f\in F}\{I^x_2(f)\}.
\end{eqnarray}
\end{Them}
\begin{Them}
Assume the conditions of Theorem 1.1.\  $\sigma $ and ${b} $ are
bounded and satisfy
  the following
\begin{eqnarray}
&&| {b}(x)-{b}(y)| + \|\sigma (x)-\sigma(y) \|+ \|(\triangledown
\sigma \cdot \sigma)(x)- (\triangledown \sigma
\cdot \sigma )(y)\|\nonumber\\
&&\|(\triangledown \sigma  \cdot \triangledown \sigma \cdot
\sigma)(x)- (\triangledown \sigma \cdot \triangledown \sigma \cdot
\sigma )(y)\| + \| (\triangledown \sigma \cdot {b})(x)-
(\triangledown \sigma \cdot {b})(y)\|\nonumber\\
&&+\| (\sigma ^T\cdot\triangledown^2 \sigma \cdot \sigma)(x)-
(\sigma^T\cdot \triangledown^2 \sigma\cdot \sigma)(y)\|\leq C |x-y|
\end{eqnarray}
for some constants $C>0 $, where $ \sigma^T$ denotes transpose of
$\sigma$,  $\triangledown \sigma $ and  $ \triangledown^2 \sigma$
denote $\sigma$'s   derivatives of first and second order with
respect to spatial variable $x$, respectively.
 Then for any random variable $X^\varepsilon_0 $ with
 ${\bf P}\{ X^\varepsilon_0\in \bar{\mathcal{O}}
 \}=1 $ and the family $\{X^\varepsilon_0, \varepsilon >0   \}  $
 satisfies for $x_0 \in \bar{\mathcal{O}}$ and any $ \delta >0  $
\begin{eqnarray}
\limsup_{\varepsilon\rightarrow 0}\varepsilon\log {\bf
P}\{|X^\varepsilon_0-x_0|>\delta  \}=-\infty ,
\end{eqnarray}
 the processes $\{ X^\varepsilon_t( X^\varepsilon_0  ):
\varepsilon >0 \} $ satisfy the large deviation principle on $E$
with
 good rate function $ I^{x_0}$ given by (1.8). In other words, for any
 open set $G$ and any closed set $F$ of $E$, we have
\begin{eqnarray}
&& \liminf_{\varepsilon \rightarrow 0}\varepsilon \log {\bf P}\{
X^\varepsilon(X^\varepsilon_0)\in G\}\geq -\inf_{f\in G}\{I^{x_0}(f)\},\\
&& \limsup_{\varepsilon \rightarrow 0}\varepsilon \log {\bf P}\{
X^\varepsilon(X^\varepsilon_0)\in F\}\leq  -\inf_{f\in
F}\{I^{x_0}(f)\}.
\end{eqnarray}
\end{Them}

This paper is organized as follows. We will study the skeleton
equation (1.5) and its Euler approximation in section 2; In section
3 we will give an exponential approximation of the solution $
\{X^\varepsilon_t \}_{t\in [0,1]} $. In sections 4 and 5, we will
prove Theorems 1.1 and 1.2.\vskip 14pt

 Throughout this paper we make the following convention: the
letter $c$ or $c(p_1, p_2, p_3, \cdots )$ depending only on $ p_1,
p_2, p_3, \cdots $ will denote an unimportant positive constant,
whose values may change from one line to another one.
\setcounter{equation}{0}
\section{{\small {\bf Skeleton Equations (1.5)}} }

\noindent In this section we mainly  study  existence, uniqueness
and approximation of solution of the skeleton equations (1.5). The
first result is the following.
\begin{Them}
Assume that $\mathcal{O} $ is a smooth bounded open set in $\Re^d$
satisfying (1.1)-(1.2) and (1.9). Let $b$ and $\sigma $ be bounded
and uniformly Lipschitz. Then for any $x\in \bar{\mathcal{O}}$  and
any $\psi \in L^2([0,1]; \Re^d ) $ the skeleton equation (1.5) has a
unique solution $ (z^\psi_t(x), k^\psi_t(x))$ and
\begin{eqnarray}
\| z^{\psi_1}-z^{\psi_2} \|^2_t\leq c(r) \int^t_0 \|\psi_1 -\psi_2
\|_s^2ds
\end{eqnarray}
holds for some constants $c(r)$, $ \psi_1, \psi_2 \in L^2( [0,1];
\Re^d ) $ with $ \| \psi_1\|^2_{L^2 }\leq r       $ and $ \|
\psi_2\|^2_{L^2 }\leq r       $,  where  $ \|f\|_t:= \sup\limits_{
s\in [0,t]}\{|f(s)|\}$ and $r>0$.
\end{Them}
\noindent {\bf Proof.}  For $\psi \in L^2([0,1]; \Re^d)    $, define
$ F(\cdot ): H:=\{ f\in E_1; \|f\|_{E_1}< +\infty \}\rightarrow
\Re^d $ by
\begin{eqnarray} \left \{\begin{array}{l}
F(z_t)=x + \int^t_0 b(z_s)ds +\int^t_0 \sigma
(z_s))\psi(s)ds-k^\psi_t(x), \\
k^\psi_t(x)=\int^t_0\xi(F(z_s))d|k^\psi(x)|_s,\\
|k^\psi(x)|_t=\int^t_0I_{\{s:F(z_s))\in \partial \mathcal{O}
\}}d|k^\psi(x)|_s,
\end{array} \right.
\end{eqnarray}
where $E_1=C([0,1]; \Re^d )$. Then there exists a constant
$c=c(\psi)>0$ such that $ \forall t\in [0,1], \forall z, z' \in H $
we have
\begin{eqnarray}
\| X-X'\|^4_t\leq c \int^t_0 \|z-z'\|^4_s ds,
\end{eqnarray}
where $ X=F(z) $ and $ X'=F(z')$.\\

We first prove (2.3).

 By chain rule,  we have for $ \phi \in C^2_b
(\bar{\mathcal{O}}) $
\begin{eqnarray}
\phi (X_t) &=&  \phi (x)  + \int^t_0 \triangledown  \phi(X_s)\cdot
b(z_s)ds +\int^t_0 \triangledown \phi(X_s) \cdot \sigma
(z_s)\psi(s)ds\nonumber\\
&&-\int^t_0\triangledown \phi(X_s) \cdot \xi(X_s)d|k^\psi(x)|_s,\\
\phi (X'_t)& =&  \phi (x)  + \int^t_0 \triangledown  \phi(X'_s)\cdot
b(z'_s)ds +\int^t_0 \triangledown \phi(X'_s) \cdot \sigma
(z'_s)\psi(s)ds\nonumber\\
&& -\int^t_0\triangledown \phi(X'_s) \cdot \xi(X'_s)d|k'^\psi(x)|_s.
\end{eqnarray}
Using conditions (1.1),(1.2) and (1.9) on $ O$,
\begin{eqnarray}
&&\frac{1}{\alpha}( \triangledown \phi (z_s ),  \xi(z_s) )|
z_s-z'_s|^2 - (z_s-z'_s, \xi(z_s) )\leq 0, \ d|k^\psi(x)|_s,
a.s., \nonumber\\
\\
&&\frac{1}{\alpha}( \triangledown \phi (z'_s ),  \xi(z'_s) )|
z_s-z'_s|^2 - (z_s-z'_s, \xi(z'_s) )\leq 0, \ d|k'^\psi(x)|_s,
a.s..\nonumber\\
\end{eqnarray}
We now calculate the following term:
$$  \exp\big \{  -\frac{2}{\alpha}[\phi(X_t) +\phi(X'_t)]\big \}|X_t-X'_t|^2.$$
Let $ f(t)= \exp\big \{  -\frac{2}{\alpha}[\phi(X_t)
+\phi(X'_t)]\big \}$. By chain rule we have
\begin{eqnarray}
f(t)|X_t&-&X'_t|^2\nonumber\\
&=&2\int^t_0f(s)(X_s-X'_s)\cdot (
b(z_s)-b(z'_s))ds\nonumber\\
&+& 2\int^t_0f(s)(X_s-X'_s)\cdot (
\sigma(z_s)-\sigma(z'_s))\cdot \psi(s)ds\nonumber\\
&-&\frac{2}{\alpha}\int^t_0f(s)|X_s-X'_s|^2\big [ \triangledown
\phi(X_s)\cdot b(z_s)  + \triangledown \phi (X'_s)\cdot
b(z'_s)\nonumber\\
& +& \triangledown \phi (X_s)\cdot \sigma(z_s) \psi (s) +
\triangledown
\phi \phi (X'_s)\cdot \sigma(z'_s)   \psi (s) \big ]ds\nonumber\\
&+&  2\int^t_0 f(s)\big [\frac{1}{\alpha}( \triangledown \phi (z_s
), \xi(z_s) )| z_s-z'_s|^2 - (z_s-z'_s, \xi(z_s) )\big
]d|k^\psi(x)|_s \nonumber\\
&+& 2\int^t_0 f(s)\big [\frac{1}{\alpha}( \triangledown \phi (z'_s
), \xi(z'_s) )| z_s-z'_s|^2 - (z_s-z'_s, \xi(z'_s) )\big
]d|k'^\psi(x)|_s.\nonumber\\
\end{eqnarray}
Since $b $ and $ \sigma $ are bounded and uniformly Lipschitz,  $
\phi \in C^2_b (\bar{\mathcal{O}}) $, we deduce from (2.6)-(2.8)
that there exists a constant $c>0$ such that
\begin{eqnarray}
|X_t-X'_t|^2&\leq & c\int^t_0 |X_s-X'_s||z_s-z'_s|[ 1+
|\psi(s)|]ds\nonumber\\
 &&+c\int^t_0 |X_s-X'_s|^2ds.
\end{eqnarray}
Consequently,
\begin{eqnarray}
|X_t-X'_t|^4&\leq & c(  1 + \|\psi \|^2_{L^2})\int^t_0
|X_s-X'_s|^4ds \nonumber\\ && + c(  1 + \|\psi \|^2_{L^2})\int^t_0
|z_s-z'_s|^4ds.\nonumber\\
\end{eqnarray}
Letting $ c_1 =c(  1 + \|\psi \|^2_{L^2})\exp\{ c(  1 + \|\psi
\|^2_{L^2})   \}$, by the Gronwall lemma, we obtain that
\begin{eqnarray}
|X_t-X'_t|^4\leq c_1 \int^t_0 |z_s-z'_s|^4ds.
\end{eqnarray}
So the proof of (2.3) has been done.\\

\noindent Next we turn to proving existence of solution of (1.5).

 Define  a
sequence $ \{ Y_n(s)\}_{n=1}^\infty$ of functions by
\begin{eqnarray} \left \{\begin{array}{l}
Y_0(t)=x, \\
Y_n(t)=F(Y_{n-1}(t)), \ t\in [0,1].
\end{array} \right.
\end{eqnarray}
By (2.3), for large enough $n$ and some constants $c>0$, we have
$$\|Y_n-Y_{n-1}\|_{E_1} \leq \frac{c}{n^2} .  $$
Define $Y(t)$ on $[0,1]$ by
$$ Y(t)=Y_0(t) + \sum^\infty_{n=1} ( Y_{n}(t)-Y_{n-1}(t)).$$
Then
$$ \|Y_n-Y\|_{E_1}\rightarrow 0\  \  \mbox{as $n \rightarrow +\infty
$}.$$ Using the first equation in (2.2) and (2.12), there exists a $
k\in  C([0,1]; \Re^d ) \cap BV([0,1])$ ( $ BV([0,1])$ denotes the
space of bounded variation functions on $[0,1]$ ) such that
$$ \|k_n-k\|_{E_1} \rightarrow 0\  \  \mbox{as $n \rightarrow +\infty
$},$$ where $k_n$ is the local time of $ Y_n$. So we can assume that
$ (Y_n, k_n) $  converges uniformly to $(Y, k) $, which are
continuous on $[0,1]$.  Letting $n\rightarrow +\infty $ in (2.12) we
know that $ (Y, k) $ solves the first equation in (1.5). The remains
of proving existence  is to check the second and third equations in
(1.5). Let $ f\in C([0,1];\Re^d )$, $ 0\leq f\leq 1$, $f|_A=1$, $A$
being any compact subset of $O$. By Fatou lemma,
$$ 0\leq \int^t_0 f(Y_s)d |k|_s \leq \liminf_{n\rightarrow
+\infty}\int^t_0 f(Y_n(s))d|k_n|_s=0,$$ which implies that
$$ \int^t_0 I_{\{s: Y_s\in \partial O\}}d|k|_s=0.$$
Thus $(Y, k ) $ solves the second equation in (1.5).

Using conditions (1.1),(1.2) and (1.9) on $O$ again, for any $ f\in
C([0,1]; \Re )$ with $f\geq 0 $ and any $ \beta \in
\bar{\mathcal{O}}$, we have
$$\int^t_0 f(s) ( Y_n(s)-\beta, dk_n(s))
+ C_0 \int^t_0 | Y_n(s)-\beta |^2 f(s)d|k_n|_s \geq 0. $$
Noticing
that the measure $ d|k_n|$ converges to some measure $da_s$ as $ n
\rightarrow +\infty $, we get that
$$\int^t_0 f(s) ( Y(s)-\beta, dk(s))
+ C_0 \int^t_0 | Y(s)-\beta |^2 f(s)da_s \geq 0. $$ Since $d|k|_s
\leq da_s$, we have $ dk_s=h_s da_s$, $h_s$ being a bounded
measurable function, thus
$$ ( Y(s)-\beta, h_s)
+ C_0| Y(s)-\beta |^2 \geq 0,$$ which, together with (1.1) and
(1.2), implies that $ h_s\in \lambda {\bf n}(Y_s) $ for some
$\lambda \geq 0 $. Hence, we find a $\xi(Y_s)\in {\bf n}(Y_s) $ such
that
$$ k_t = \int^t_0\xi(Y_s)d|k|_s,$$
that is, $ (Y,k) $ also solves the third equation in (1.5).
Therefore we have proved the existence of solution of (1.5).\\

 Let $ (Y, k) $ and $(Y', k' ) $ be two solutions of (1.5) with
 $Y(0)=Y'(0)=x $ and $k_0=k'_0 =0 $. By (2.3),
 $$ \|Y-Y'\|_{E_1}^4 \leq c \int^1_0 \| Y-Y'\|^4_s ds.$$
So $ Y=Y'$. We then deduce from the first equation in (1.5) that $
k=k'$. This implies  the uniqueness of solution of Eq.(1.5).\\

To the proof of Theorem 2.1 end, we  to prove (2.1). Let $ \psi_i
\in L^2([0,1]; \Re^d) $ and $ (z^{\psi_i}_t(x), k^{\psi_i}_t(x))$ be
the solutions of Eq.(1.5) corresponding to $\psi_i $ for $i=1,2$.
Replacing $ f(t) $ by $  \exp\big \{
-\frac{2}{\alpha}[\psi_1(z^{\psi_1}_t(x)) + \psi_2(z^{\psi_2}_t(x))
]\big \}$ in (2.8) and  using the same way as in (2.10), we can
prove that for any $r>0$ and any $ \| \psi_i\|^2_{L^2 }\leq r $  for
$i=1,2,$
\begin{eqnarray}
| z^{\psi_1}_t(x)-z^{\psi_2}_t(x)|^2&\leq & c(r)\int^t_0 |
z^{\psi_1}_s(x)-z^{\psi_2}_s(x) |^2\nonumber\\&& +c(r) \int^t_0 |
\psi_1(s)-\psi_2(s)|^2ds.
\end{eqnarray}
 Thus by  Gronwall lemma  we complete the proof. \qquad
$\Box$\vskip 19pt

For $n\geq 1$, let $ \phi_n   $ be a step function on $[0,1]$
defined by $ \phi_n(t)=\frac{k}{2^n} $ for $t\in [ \frac{k}{2^n},
\frac{k+1}{2^n}]$, $ k=0,1,2,\cdots, 2^n-1.$  Let $\{ z_n^\psi  \}
\in C([0,1]\times \bar{\mathcal{O}}; \bar{\mathcal{O}}    ) $ be the
Euler approximating sequence of Eq.(1.5) defined by the following
\begin{eqnarray} \left \{\begin{array}{l}
z^\psi_n(t,x)=x + \int^t_0 b(z^\psi_n(\phi_n(s),x))ds +\int^t_0
\sigma(z^\psi_n( \phi_n(s), x))\psi(s)ds-k^\psi_n(t), \\
k^\psi_n(t)=\int^t_0\xi(z^\psi_n(s,x))d|k^\psi_n|_s,\\
|k^\psi_n|_t=\int^t_0I_{\{s:z^\psi_n(s,x)\in \partial \mathcal{O}
\}}d|k^\psi_n|_s
\end{array} \right.
\end{eqnarray}
for $n\geq 1$ and $\psi\in L^2([0,1]; \Re^d ).  $  The second result
of this section is the following.
\begin{Them} Let $O$ satisfy conditions (1.1)-(1.3). Then for any $r>0$ we have
\begin{eqnarray}
\lim_{n\rightarrow \infty}\sup_{\{  \psi: \ I_1(\psi)\leq r \}}\|
z_n^\psi-z^\psi\|_{E_2}=0,
\end{eqnarray}
where $ \| \cdot \|_{E_2} $ denotes the uniform norm of $E_2:=
C([0,1]\times \bar{\mathcal{O}}; \bar{\mathcal{O}}    )     $.
\end{Them}
{\bf Proof.} By (1.5) and (2.14)
\begin{eqnarray}
z_n^\psi(t,x)-z^\psi_t(x)&=&\int^t_0 \big [ b(z^\psi_n(\phi_n(s),x))
-b(z^\psi_s(x))\big ]ds \nonumber\\
&+&\int^t_0 \big [ \sigma(z^\psi_n(\phi_n(s),x))
-\sigma(z^\psi_s(x))\big ]\psi(s)ds \nonumber\\
&-& \big [ k^\psi_n(t)-k^\psi(s)     \big]:= \omega (t)-L^\psi_n(t).
\end{eqnarray}
Since $I_1( \psi )\leq r   $, $b$ and $ \sigma $ are bounded, we
know  that $ \omega \in BV([0,1])$. Using Theorem 2.1 in \cite{s3} ,
we get that
\begin{eqnarray}
|L^\psi_n|_t &\leq & \int^t_0 |\big [ b(z^\psi_n(\phi_n(s),x))
-b(z^\psi_s(x))\big ]|ds \nonumber\\
&+&\int^t_0| \big [ \sigma(z^\psi_n(\phi_n(s),x))
-\sigma(z^\psi_s(x))\big ]|| \psi(s)|ds.
\end{eqnarray}
So we deduce from (2.16) that
\begin{eqnarray}
| z_n^\psi(t,x)-z^\psi_t(x)|&\leq & 2\int^t_0 |\big [
b(z^\psi_n(\phi_n(s),x))
-b(z^\psi_s(x))\big ]|ds \nonumber\\
&+&2\int^t_0| \big [ \sigma(z^\psi_n(\phi_n(s),x))
-\sigma(z^\psi_s(x))\big ]|| \psi(s)|ds.\nonumber\\
\end{eqnarray}
Similarly,
\begin{eqnarray}
| z_n^\psi(t,x)-z^\psi_n(\phi_n(t),x
 )|&\leq & 2\int^t_{\phi_n(t) }  |b(z^\psi_n(\phi_n(s),x))
ds\nonumber\\
&& +2\int^t_{\phi_n(t)}|  \sigma(z^\psi_n(\phi_n(s),x))
|| \psi(s)|ds.\nonumber\\
&\leq & \frac{2c}{2^n}\big [ 1 + | \int^1_0 | \psi(s)   | \big
]ds\nonumber\\
&\leq & c(r)\frac{1}{2^n}.
\end{eqnarray}
Hence,  using (2.18) and (2.19) as well as $b$ and $\sigma $ are
uniformly Lipschitz  we have
\begin{eqnarray}
&&\sup_{\{x\in\bar{\mathcal{O}}, s\leq t        \} }\{|
z_n^\psi(s,x)-z^\psi_s(x
 )|\}\nonumber \\
 &\leq & c \int^t_0 [1+ | \psi(s) |][\sup_{\{x\in\bar{\mathcal{O}}, u\leq s     \}    }\{|
z_n^\psi(u,x)-z^\psi_u(x
 )|\}]ds\nonumber\\
 &&+ c(r)\frac{1}{2^n}\int^t_0 [1+ | \psi(s) |]ds,
\end{eqnarray}
which, together with Gronwall lemma, implies that
\begin{eqnarray}
\| z_n^\psi-z^\psi\|_{E_2}\leq c(r)\frac{1}{2^n}\int^1_0 [1+ |
\psi(s) |]ds \exp \{c\int^1_0 [1+ | \psi(s) |]ds \}.\nonumber\\
\end{eqnarray}
Since $\int^1_0 [1+ | \psi(s) |]ds \leq ( 1+ \sqrt{2r}       ) $,
letting $ n \rightarrow  \infty $ in (2.21) we have
\begin{eqnarray*}
\lim_{n\rightarrow \infty}\sup_{\{  \psi: \ I_1(\psi)\leq r \}}\|
z_n^\psi-z^\psi\|_{E_2}=0.
\end{eqnarray*}
Thus we complete the proof of Theorem 2.2. \quad $\Box$\vskip 19pt

We define $F_n(\cdot ) $ from $E_1$ to $E_2$ by
\begin{eqnarray} \left \{\begin{array}{l}
F_n(\omega )(0,x)=x ,\\ F_n( \omega  )(t,x)= F_n(  \omega
)(\frac{k}{2^n},  x) + b(F_n(\omega)(\frac{k}{2^n},
x))(t-\frac{1}{2^n}) \\ \qquad  \qquad  \quad +
\sigma(F_n(\omega)(\frac{k}{2^n}, x))(\omega (t)- \omega
(\frac{1}{2^n})) \\ \qquad  \qquad  \quad -
\int^t_{\frac{k}{2^n}}\xi( F_n( \omega  )(s,x)     ) d |L^x_n|_s,\\
L^x_n(s)=\int^t_0\xi( F_n( \omega  )(s,x)     ) d |L^x_n|_s,\\
|L^x_n|_t =\int^t_0I_{\{s: F_n( \omega  )(s,x) \in \partial
\mathcal{O} \}}d|L^x_n|_s
\end{array} \right.
\end{eqnarray}
for $n\geq 1$ and $k=0,1, \cdots, 2^n-1$.
 The third result of this
section is the following.
\begin{Them}
Let $O$ satisfy conditions (1.1)-(1.3) and (1.9). Then for $n\geq 1$
the maps  $F_n(\cdot ) $ from $E_1$ to $E_2$ are continuous.
\end{Them}
Before proving proof of Theorem 2.3, we need the following lemmas.
\begin{Lemma}
Let the smooth bounded $O$ satisfy (1.1)-(1.3) and $ w \in E.$
Assume that for every $ x\in \bar{\mathcal{O}}$, \  $( Y_x, k_x ) $
is the unique solution of the following Skorohod equation(
abbreviated by $(\omega_x, O, {\bf n})$),
\begin{eqnarray} \left \{\begin{array}{l}
Y_x(t)=\omega (t,x)-k_x(t) ,\\
k_x(t)=\int^t_0\xi( Y_x(s)) d |k_x|_s,\\
|k_x|_t =\int^t_0I_{\{s: Y_x(s) \in \partial \mathcal{O}
\}}d|k_x|_s.
\end{array} \right.
\end{eqnarray}
Then
\begin{eqnarray}\|k\|_{E_2}\leq M(\alpha, C_0,
\bar{\mathcal{O}}) < +\infty,
\end{eqnarray}
 where $M(\alpha, C_0, \bar{\mathcal{O}})$ are some positive
constants depending only on $ \alpha $, $ C_0 $ and $
\bar{\mathcal{O}}.    $
\end{Lemma}
{\bf Proof.} Since $\bar{\mathcal{O}} $ satisfies (1.3) and is
bounded, the proof is completely similar to  that of Lemma 1.2 in
\cite{s3}, we omit it here. \quad \quad $\Box $
\begin{Lemma}
Let the smooth bounded $O$ satisfy (1.1) and (1.2). Assume that
$(Y_i, k_i)$ , $i=1,2$, are unique solutions of Skorohod equations $
( w_1, O, {\bf n}) $ and $ ( w_2, O, {\bf n})$, respectively.  Then
\begin{eqnarray}
\| Y_1-Y_2\|_t^2 &\leq &\big [1+ \exp\{2C_0 (|k_1|_t +|k_2|_t)) \}
\big ]\|\omega_1-\omega_2\|^2_t\nonumber\\
&& + \frac{3}{C_0}\exp\{4C_0 (|k_1|_t +|k_2|_t))
\}\|\omega_1-\omega_2\|_t.
\end{eqnarray}
\end{Lemma}
{\bf Proof. } Since the first inequality (6) in Lemma 1.1 in
\cite{s3} holds for $\omega \in E $ due to constant $C$ there can be
replaced by $2C_0$, we have
\begin{eqnarray}
|Y_1(t)-Y_2(t)|^2 &\leq & |\omega_1(t)-\omega_2(t)|^2 +
2(k_1(t)-k_2(t))\cdot (\omega_1(t)-\omega_2(t) )\nonumber\\
&& + \exp\big \{2C_0(|k_1|_t +|k_2|_t )\big \}\cdot \big \{
\|\omega_1-\omega_2\|^2_t
\nonumber\\
&&+\frac{1}{C_0}\|\omega_1-\omega_2\|_t +\|(\omega_1-\omega_2 )\cdot
(k_1-k_2)\|_t \big \},
\end{eqnarray}
which implies that (2.25). \quad \quad $ \Box $ \\

Now we turn to proving Theorem 2.3.\\

\noindent {\bf Proof of Theorem 2.3. } For $ \omega_i \in E_1$ with
$ \|\omega_i \|_{E_1} < + \infty $, $i=1,2$, define $( F_n(
\omega_i),L^{x,i}_n) $ by (2.22). Let
\begin{eqnarray}
\widetilde{w_i}(t,x)&=&F_n(  \omega_i )(\frac{k}{2^n},  x) +
b(F_n(\omega_i)(\frac{k}{2^n}, x))(t-\frac{1}{2^n}) \nonumber\\
&& + \sigma(F_n(\omega_i)(\frac{k}{2^n}, x))(\omega_i (t)- \omega_i
(\frac{1}{2^n})).
\end{eqnarray}
Then
\begin{eqnarray} \left \{\begin{array}{l}
 F_n( \omega_i  )(t,x)= \widetilde{w_i}(t,x) - L^{x,i}_n,
\\
L^{x,i}_n(s)=\int^t_0\xi( F_n( \omega_i  )(s,x)     ) d |L^{x, i}_n|_s,\\
|L^{x, i}_n|_t =\int^t_0I_{\{s: F_n( \omega_i  )(s,x) \in \partial
\mathcal{O} \}}d|L^{x,i}_n|_s.
\end{array} \right.
\end{eqnarray}
By Lemma 2.1, $ \exists \  M_1=M_1(  \alpha, O, C_0) <+ \infty     $
such that for $i=1,2$,
\begin{eqnarray}
\sup_{n\geq 1}\big \{\| L^{\cdot, i}_n\|_{E_2}\big \} \leq M_1.
\end{eqnarray}
Setting $ Y_n(t, x)= F_n(\omega_1)(t,x)-F_n( \omega_2)(t,x)$.  For
$t\in [0, \frac{1}{2^n}]$, by Lemma 2.2 and (2.4), $\exists M_2 =
M_2( \alpha, O, C_0 )>0 $  such that
\begin{eqnarray}
\| Y_n(\cdot, x)\|_{\frac{1}{2^n}}\leq M_2\big\{ \|
\widetilde{w_1}(\cdot , x)-\widetilde{w_2}(\cdot, x
)\|^2_{\frac{1}{2}}+ \|\widetilde{w_1}(\cdot ,
x)-\widetilde{w_2}(\cdot, x )       \|_{\frac{1}{2} }
\}^{\frac{1}{2} },\nonumber\\
\end{eqnarray}
where $ \|f\|_t =\sup_{s\leq t}\{ |f(s)|  \}$. Noticing that for $
t\in [0, \frac{1}{2^n}]$,
$$ \widetilde{w_1}(\cdot , x)-\widetilde{w_2}(\cdot, x
) = \sigma (x)[
(\omega_1(t)-\omega_2(t))-(\omega_1(0)-\omega_2(0))],
$$ we have
\begin{eqnarray}
\| \widetilde{w_1}(\cdot , x)-\widetilde{w_2}(\cdot, x )
\|_{\frac{1}{2^n}}\leq c \| \omega_1-\omega_2 \|_{\frac{1}{2^n}}.
\end{eqnarray}
Using (2.30) and (2.31), we have
\begin{eqnarray}
\| Y_n(\cdot, x)\|_{\frac{1}{2^n}}\leq c( \alpha, O, C_0) \big\{ \|
\omega_1-\omega_2\|^2_{\frac{1}{2}}+ \|
\omega_1-\omega_2\|_{\frac{1}{2} } \}^{\frac{1}{2} }.
\end{eqnarray}
Since
\begin{eqnarray}
\| Y_n(\cdot, x)\|_{\frac{2}{2^n}} &\leq &  \| Y_n(\cdot,
x)\|_{\frac{1}{2^n}}\nonumber\\
&+&\sup_{t\in \bigtriangleup_2}\{| \int^t_{\frac{1}{2^n}} \big
[b(F_n(\omega_1)(\frac{k}{2^n}, x))
 -b(F_n(\omega_2)(\frac{k}{2^n}, x))      \big ]ds   |\}\nonumber\\
 &+&\sup_{t\in \bigtriangleup_2}\{|  \big
[\sigma(F_n(\omega_1)(\frac{k}{2^n}, x))
 -\sigma (F_n(\omega_2)(\frac{k}{2^n}, x))      \big ](\omega_1(t) - \omega_1(0)     )
 |\}\nonumber\\
 &+&\sup_{t\in \bigtriangleup_2}\{| \sigma(F_n(\omega_2)(\frac{k}{2^n},
x))\big [ (\omega_1(t)-  \omega_2(t))-(  \omega_1(0)- \omega_2(0))
\big ] |\},\nonumber\\
\end{eqnarray}
where $ \bigtriangleup_2:= [ \frac{1}{2^n}, \frac{2}{2^n}] $, we
have
\begin{eqnarray*}
\| Y_n(\cdot, x)\|_{\frac{2}{2^n}}\leq c( \alpha, O, C_0, \|\omega_1
\|_{E_1})\| \omega_1  - \omega_2   \|^{\frac{1}{2}}_{E_1}( 1+ \|
\omega_1 - \omega_2   \|^{\frac{1}{2}}_{E_1}   ) +c\| \omega_1  -
\omega_2 \|_{E_1}.
\end{eqnarray*}
Doing the same procedure as in estimating $\| Y_n(\cdot,
x)\|_{\frac{2}{2^n}} $, we can find $c_1=c_1( \alpha, O, C_0,
\|\omega_1 \|_{E_1})$ and $c_2=c_2( \alpha, O, C_0, \|\omega_1
\|_{E_1})$
such that
\begin{eqnarray*}  \|Y_n(\cdot, \cdot )\|_{E_2} \leq c_1  \| \omega_1  - \omega_2
\|^{\frac{1}{2}}_{E_1}( 1+ \| \omega_1  - \omega_2
\|^{\frac{1}{2}}_{E_1} ) +c_2\| \omega_1  - \omega_2 \|_{E_1}.
\end{eqnarray*}
Thus we complete the proof by the last inequality. \quad $\Box
$\vskip 19pt
 \setcounter{equation}{0}
\section{{\small {\bf Exponential approximation of the solution  for SDE(1.4)}} }
 \noindent We consider the following Euler approximation of SDE
(1.4),
\begin{eqnarray} \left \{\begin{array}{l}
X^\varepsilon_n(t,x)=x + \int^t_0 b(X^\varepsilon_n(\phi_n(s),x))ds
+\sqrt{\varepsilon}\int^t_0 \sigma
(X^\varepsilon_n(\phi_n(s),x))dB_s\\
\qquad \qquad
-L^\varepsilon_n(t,x), \\
L^\varepsilon_n(t,x)=\int^t_0\xi(X^\varepsilon_n(s,x))d
|L^\varepsilon_n(x)|_s,\\
|L^\varepsilon_n(x)|_t=\int^t_0I_{\{s:X^\varepsilon_n(s,x)\in
\partial \mathcal{O} \}}d|L^\varepsilon_n(x)|_s
\end{array} \right.
\end{eqnarray}
for $ n\geq 1$ and $x\in \bar{\mathcal{O}    }$. The main result of
this section is the following.
\begin{Them}
Let $O$ satisfy conditions(1.1)-(1.3) and (1.9). Then  we have
\begin{eqnarray}
\lim_{n\rightarrow \infty }\limsup_{ \varepsilon\rightarrow 0  }
\varepsilon \log {\bf P} \big\{ \sup_{t\in [0,1]}|
X^\varepsilon_t(x)-X^\varepsilon_n(t,x)|\geq \delta \big \}=-\infty
\end{eqnarray}
for any $ \delta>0  . $
\end{Them}
{\bf Proof.} We first define stopping time  $ \tau_1 $ for
$\delta_1>0$ by
$$ \tau_1= \inf\big \{t\geq 0: |X^\varepsilon_n(t)-
X^\varepsilon_n(\phi_n(t))|\geq \delta_1\big\}\wedge 1.           $$
Then
\begin{eqnarray}
{\bf P}\{ \tau_1\leq 1  \}\leq \sum^{2^n-1}_{k=0}{\bf P} \big \{
\sup_{t\in [\frac{k}{2^n}, \frac{k+1}{2^n}]}|X^\varepsilon_n(t)-
X^\varepsilon_n(\phi_n(t))|\geq \delta_1 \big \}.
\end{eqnarray}
Since
\begin{eqnarray}
X^\varepsilon_n(t)- X^\varepsilon_n(\phi_n(t))&=&\int^t_{\phi_n(t)}
b(X^\varepsilon_n(\phi_n(s),x))ds +\sqrt{\varepsilon}\int^t_{
\phi_n(t)}\sigma (X^\varepsilon_n(\phi_n(s),x))dB_s\nonumber\\
&&-\int^t_{\phi_n(t)}\xi(X^\varepsilon_n(s,x))d
|L^\varepsilon_n(x)|_s,
\end{eqnarray}
by using (2.24) and (2.25), there exists a  positive constant $c_1$
depending only on $\alpha$, $O$ and $C_0 $ such that
\begin{eqnarray}
|X^\varepsilon_n(t)&-& X^\varepsilon_n(\phi_n(t))|\nonumber\\
&\leq  &  c \big [ | \int^t_{\phi_n(t)}
b(X^\varepsilon_n(\phi_n(s),x))ds +\sqrt{\varepsilon}\int^t_{
\phi_n(t)}\sigma (X^\varepsilon_n(\phi_n(s),x))dB_s |\nonumber\\
&+& | \int^t_{\phi_n(t)} b(X^\varepsilon_n(\phi_n(s),x))ds
+\sqrt{\varepsilon}\int^t_{ \phi_n(t)}\sigma
(X^\varepsilon_n(\phi_n(s),x))dB_s  |^{\frac{1}{2}} \big
]\nonumber\\
&\leq & \big [ \frac{c_2}{2^n} + c_2\sqrt{\varepsilon} \max_{0\leq
t\leq \frac{1}{2^n}}\{ | \widetilde{B}_t|\}+\sqrt{\frac{c_2}{2^n} +
c_2\sqrt{\varepsilon} \max_{0\leq t\leq \frac{1}{2^n}}\{ |
\widetilde{B}_t|\}}\big ]\nonumber\\
&:= & b_n +\sqrt{b_n},
\end{eqnarray}
where $ \widetilde{B}_t=B(t+\frac{k}{2^n} -B( \frac{k}{2^n} )   $ is
 also  a Brownian motion. So, by choosing $\frac{\delta_1}{2} \leq 1
 $, we have
\begin{eqnarray}
&&{\bf P} \big \{ \sup_{t\in [\frac{k}{2^n},
\frac{k+1}{2^n}]}|X^\varepsilon_n(t)-
X^\varepsilon_n(\phi_n(t))|\geq \delta_1 \big \}\nonumber\\
&& \leq {\bf P}\big \{ b_n\geq \frac{\delta_1}{2}  \big \} + {\bf
P}\big \{ b_n\geq (\frac{\delta_1}{2})^2  \big \}\nonumber\\
&& \leq 2 {\bf P}\big \{ b_n\geq (\frac{\delta_1}{2})^2  \big
\}\nonumber\\
&& \leq 2 {\bf P}\big \{ \max_{0\leq t\leq \frac{1}{2^n}}\{ |
\widetilde{B}_t|\}\geq  \frac{ (\frac{\delta_1^2}{4}-
\frac{c_1}{2^n} ) }{\sqrt{\varepsilon}c_1 } \big\}\nonumber\\
&& \leq 8d \exp\big\{ -\frac{2^n  (\frac{\delta_1^2}{4}-
\frac{c_1}{2^n} )^2 }{2d\varepsilon c_1^2}         \big \}
\end{eqnarray}
by Lemma 5.2.1 in \cite{s10}. So
\begin{eqnarray*}
{\bf P}\{ \tau_1\leq 1  \}\leq 2^n 8d \exp\big\{ -\frac{2^n
(\frac{\delta_1^2}{4}- \frac{c_1}{2^n} )^2 }{2d\varepsilon c_1^2}
\big \}.
\end{eqnarray*}
Consequently,
\begin{eqnarray*}
\varepsilon \log{\bf P}\{ \tau_1\leq 1  \}\leq  \varepsilon n \log 2
+ \varepsilon\log ( 8d) - 2^n \frac{ (\frac{\delta_1^2}{4}-
\frac{c_1}{2^n} )^2 }{2d c_1^2}.
\end{eqnarray*}
Hence,
\begin{eqnarray}
\lim_{n\rightarrow \infty }\limsup_{\varepsilon\rightarrow
0}\varepsilon \log{\bf P}\{ \tau_1\leq 1 \}=-\infty .
\end{eqnarray}
Next we define functions $ \Psi $  and $\Phi_\lambda  $ on $[0,
+\infty )$ by
$$ \Psi (x)= \int^x_0 \frac{ds}{\delta_1^2+ s} \ \mbox{ and } \
\Phi_\lambda(x)=\exp\big \{ \lambda \Psi(x)  \big \} \ \mbox{for
$\lambda >0$ .}
$$
For $\delta >0$, define stopping time $ \tau_2$ by
 $$ \tau_2=\inf\{
t>0: |z_t| \geq \delta \}\wedge \tau_1  ,$$
where
$$ z_t= X^\varepsilon_t(x)-X^\varepsilon_n(t,x). $$
Let
\begin{eqnarray}
m_t&=&z_{t\wedge \tau_2 },\ f(x)=|x|^2, \ x\in \Re^d,\\
\widetilde{b}_t&=&b (X^\varepsilon_{t\wedge \tau_2}(x)  )
-b(X^\varepsilon_n(\phi_n(t),x)    ),\\
\widetilde{\sigma}_t&=&\sqrt{\varepsilon}\big [\sigma
(X^\varepsilon_{t\wedge \tau_2}(x) )
-\sigma(X^\varepsilon_n(\phi_n(t),x)    )\big ].
\end{eqnarray}
Then
\begin{eqnarray}
|\widetilde{b}_t|&\leq & c [ \delta_1 +|m_t|],\\
\|\widetilde{\sigma}_t\|^2&\leq & c\varepsilon [\delta_1^2 +
|m_t|^2],\\
 m_t&=& \int^{t\wedge \tau_2}_0\widetilde{b}_s ds  + \int^{t\wedge
\tau_2}_0\widetilde{\sigma}_sdB_s\nonumber\\
& -& \int^{t\wedge \tau_2}_0 \xi(X^\varepsilon_{s\wedge \tau_2}(x))d
|L^\varepsilon(x)|_s\nonumber\\
& + & \int^{t\wedge \tau_2}_0\xi(X^\varepsilon_n(s\wedge\tau_2,x))d
|L^\varepsilon_n(x)|_s.
\end{eqnarray}
Define
\begin{eqnarray*}
D_t&=&\phi( X^\varepsilon_{t\wedge \tau_2}(x)   ) + \phi(
X^\varepsilon_n(t\wedge \tau_2, x)),\\
N_t&=&\exp\{-\frac{2}{\alpha}D_t \}.
\end{eqnarray*}
By It\^{o}'s formula
\begin{eqnarray}
\Phi_\lambda (f(m_t)  )&=&1+ \int^{t\wedge \tau_2}_0\Phi_\lambda'
(f(m_s)  )df(m_s)\nonumber\\
& +& \frac{1}{2}\int^{t\wedge \tau_2}_0\Phi_\lambda'' (f(m_s)
)df(m_s)\cdot df(m_s),\\
f(m_t)N_t&=&\int^{t\wedge \tau_2}_0N_sdf(m_s) + \int^{t\wedge
\tau_2}_0f(m_s)dN_s \nonumber\\
&+& \int^{t\wedge \tau_2}_0df(m_s)\cdot dN_s,
\end{eqnarray}
\begin{eqnarray}
f(m_t)&=& 2\int^{t\wedge \tau_2}_0 m_s\cdot \widetilde{b}_s ds +2
\int^{t\wedge \tau_2}_0 m_s\cdot \widetilde{\sigma}_sdB_s\nonumber\\
&-& 2\int^{t\wedge \tau_2}_0 m_s\cdot \xi(X^\varepsilon_{s\wedge
\tau_2}(x))d |L^\varepsilon(x)|_s\nonumber\\
& +&2\int^{t\wedge \tau_2}_0 m_s\cdot
\xi(X^\varepsilon_n(s\wedge\tau_2,x))d
|L^\varepsilon_n(x)|_s\nonumber\\
&+& \int^{t\wedge \tau_2}_0{\bf trace
}\{\widetilde{\sigma}_s\widetilde{\sigma}_s^T\}ds,
\end{eqnarray}
\begin{eqnarray}
 dN_t&=& -\frac{2}{\alpha}N_t dD_t +
\frac{2}{\alpha^2}N_tdD_t\cdot dD_t,\\
 dD_t &=& \big [
(\triangledown \phi^T \cdot {b} )(X^\varepsilon_{t\wedge \tau_2}(x)
)+ (\triangledown \phi^T \cdot {b} )(
X^\varepsilon_n(t\wedge\tau_2,x) )         \big ]
dt\nonumber\\
&&-\big [(\triangledown \phi^T \cdot {\xi} )(X^\varepsilon_{t\wedge
\tau_2}(x) ) d |L^\varepsilon(x)|_t \nonumber\\
&& + (\triangledown \phi^T \cdot {b} )(
X^\varepsilon_n(t\wedge\tau_2,x) )d |L^\varepsilon_n(x)|_t\big
]\nonumber\\
 &&+ \sqrt{\varepsilon}\big [  (\triangledown \phi^T \cdot {\sigma}
)(X^\varepsilon_{t\wedge \tau_2}(x) )+( \triangledown \phi^T \cdot
{\sigma} )( X^\varepsilon_n(t\wedge\tau_2,x) ) \big
]dB_t\nonumber\\
&&+\frac{\varepsilon}{2}{\bf trace}\big \{(\triangledown^2\phi \cdot
\sigma \cdot \sigma^T )( X^\varepsilon_{t\wedge \tau_2}(x)
)\nonumber\\
&& +(\triangledown^2 \phi\cdot \sigma \cdot \sigma^T
)(X^\varepsilon_n(t\wedge\tau_2,x)
    ) \big\}dt.
\end{eqnarray}
Therefore, we have the following  stochastic contractions:
\begin{eqnarray}
&& dD_t\cdot dD_t = \varepsilon {\bf trace}\{\big ( \triangledown
\phi^T \cdot {\sigma} \cdot \triangledown \phi \cdot {\sigma^T} \big
)( X^\varepsilon_{t\wedge \tau_2}(x)    )\nonumber\\
 &&\quad \qquad \qquad   \big ( \triangledown \phi^T \cdot {\sigma} \cdot \triangledown
\phi \cdot {\sigma^T} \big )( X^\varepsilon_n(t\wedge\tau_2,x)    )
\}dt,\\
&& df(m_t)\cdot df(m_t)=4{\bf trace}\{ m_t\cdot
\widetilde{\sigma_t}\cdot \widetilde{\sigma_t}^T m_t^T\}dt.
\end{eqnarray}
By using (3.17)-(3.19),
\begin{eqnarray}
dN_t &=& -\frac{2}{\alpha}N_t \sqrt{\varepsilon}\big [
(\triangledown \phi^T \cdot {\sigma} )(X^\varepsilon_{t\wedge
\tau_2}(x) )+( \triangledown \phi^T \cdot {\sigma} )(
X^\varepsilon_n(t\wedge\tau_2,x) ) \big
]dB_t\nonumber\\
&&+\frac{2}{\alpha}N_t \big [(\triangledown \phi^T \cdot {\xi}
)(X^\varepsilon_{t\wedge
\tau_2}(x) ) d |L^\varepsilon(x)|_t \nonumber\\
&& + (\triangledown \phi^T \cdot {\xi} )(
X^\varepsilon_n(t\wedge\tau_2,x) )d |L^\varepsilon_n(x)|_t\big
]\nonumber\\
&&-\frac{2}{\alpha}N_t \big [ (\triangledown \phi^T \cdot {b}
)(X^\varepsilon_{t\wedge \tau_2}(x) )+ (\triangledown \phi^T \cdot
{b} )( X^\varepsilon_n(t\wedge\tau_2,x) )         \big ]
dt\nonumber\\
&&- \frac{\varepsilon}{\alpha}N_t{\bf trace}\big
\{(\triangledown^2\phi \cdot \sigma \cdot \sigma^T )(
X^\varepsilon_{t\wedge \tau_2}(x)
)\nonumber\\
&& +(\triangledown^2 \phi\cdot \sigma \cdot \sigma^T
)(X^\varepsilon_n(t\wedge\tau_2,x)
    ) \big\}dt\nonumber\\
&& +\frac{2}{\alpha^2}N_t \varepsilon {\bf trace}\{\big (
\triangledown \phi^T \cdot {\sigma} \cdot \triangledown \phi \cdot
{\sigma^T} \big
)( X^\varepsilon_{t\wedge \tau_2}(x)    )\nonumber\\
&&+ \big ( \triangledown \phi^T \cdot {\sigma} \cdot \triangledown
\phi \cdot {\sigma^T} \big )( X^\varepsilon_n(t\wedge\tau_2,x)    )
\}dt.
\end{eqnarray}
So the stochastic contraction  $df(m_t)\cdot dN_t $ can be given by
\begin{eqnarray}
df(m_t)\cdot dN_t&=&-\frac{4}{\alpha}N_t\sqrt{\varepsilon}{\bf
trace}\big \{\widetilde{\sigma_t}^Tm_t^T \big [ (\triangledown
\phi^T \cdot {\sigma} )(X^\varepsilon_{t\wedge \tau_2}(x)
)\nonumber\\
&&+( \triangledown \phi^T \cdot {\sigma} )(
X^\varepsilon_n(t\wedge\tau_2,x) ) \big ]\big \}dt.
\end{eqnarray}
Using (3.15), (3.16), (3.21) and (3.22) we obtain that
\begin{eqnarray}
f(m_t)N_t&=&2\int^{t\wedge \tau_2}_0N_s\big\{m_s\cdot
\widetilde{\sigma_s} - \frac{\sqrt{\varepsilon}}{\alpha} f(m_s)
 \big [
(\triangledown \phi^T \cdot {\sigma} )(X^\varepsilon_{s\wedge
\tau_2}(x) )\nonumber\\&&+( \triangledown \phi^T \cdot {\sigma} )(
X^\varepsilon_n(s\wedge\tau_2,x) ) \big ]          \big
\}dB_s\nonumber\\
&&+ 2\int^{t\wedge \tau_2}_0N_s\big [\frac{1}{\alpha}|m_s|^2 \big (
\triangledown \phi( X^\varepsilon_{s\wedge \tau_2}(x)    ), {\xi}
(X^\varepsilon_{s\wedge \tau_2}(x) ) \big )\nonumber\\
&&\qquad \qquad- \big ( m_s, {\xi} (X^\varepsilon_{s\wedge
\tau_2}(x) ) \big ) \big ]d
|L^\varepsilon(x)|_s\nonumber\\
&&+ 2\int^{t\wedge \tau_2}_0N_s\big [\frac{1}{\alpha}|m_s|^2 \big (
\triangledown \phi( X^\varepsilon_n(s\wedge\tau_2,x)    ), {\xi}
(X^\varepsilon_n(s\wedge\tau_2,x) ) \big )\nonumber\\
&&\qquad \qquad- \big ( m_s, {\xi} (X^\varepsilon_n(s\wedge\tau_2,x)
) \big ) \big ]d|L^\varepsilon_n
(x)|_s\nonumber\\
&& +\int^{t\wedge \tau_2}_0N_s\big [2m_s\cdot \widetilde{b}_s + {\bf
trace }\{\widetilde{\sigma}_s\widetilde{\sigma}_s^T\} \big ]ds\nonumber\\
&&-\int^{t\wedge \tau_2}_0\frac{2}{\alpha}N_sf(m_s) \big [
(\triangledown \phi^T \cdot {b} )(X^\varepsilon_{s\wedge \tau_2}(x)
)\nonumber\\
&&\qquad   \qquad \quad + (\triangledown \phi^T \cdot {b} )(
X^\varepsilon_n(s\wedge\tau_2,x) )         \big ]
ds\nonumber\\
&&-\int^{t\wedge \tau_2}_0\frac{\varepsilon}{\alpha}N_sf(m_s) {\bf
trace}\big \{(\triangledown^2\phi \cdot \sigma \cdot \sigma^T )(
X^\varepsilon_{s\wedge \tau_2}(x)
)\nonumber\\
&& \qquad \qquad \quad +(\triangledown^2 \phi\cdot \sigma \cdot
\sigma^T )(X^\varepsilon_n(s\wedge\tau_2,x)
    ) \big\}ds\nonumber\\
&& +\int^{t\wedge \tau_2}_0 \frac{2\varepsilon}{\alpha^2}N_s  {\bf
trace}\{\big ( \triangledown \phi^T \cdot {\sigma} \cdot
\triangledown \phi \cdot {\sigma^T} \big
)( X^\varepsilon_{s\wedge \tau_2}(x)    )\nonumber\\
&&\qquad \qquad \quad + \big ( \triangledown \phi^T \cdot {\sigma}
\cdot \triangledown \phi \cdot {\sigma^T} \big )(
X^\varepsilon_n(s\wedge\tau_2,x)    )
\}ds,\nonumber\\
&& -\int^{t\wedge \tau_2}_0
\frac{4}{\alpha}N_s\sqrt{\varepsilon}{\bf trace}\big
\{\widetilde{\sigma_s}^Tm_s^T \big [ (\triangledown \phi^T \cdot
{\sigma} )(X^\varepsilon_{s\wedge \tau_2}(x)
)\nonumber\\
&&\qquad \qquad \quad +( \triangledown \phi^T \cdot {\sigma} )(
X^\varepsilon_n(s\wedge\tau_2,x) ) \big ]\big \}ds\nonumber\\
&:=& \sum^8_{i=1}a_i(t).
\end{eqnarray}
Using conditions(1.1)-(1.3) and (1.9), and $\phi \in C_b( \Re^d ) $,
$ \exists $ a positive constant $c$  such that $ \frac{1}{c}\leq N_t
\leq c$, since $ a_2(t) + a_3(t) \leq 0 $, by (3.23),  we have
\begin{eqnarray}
f(m_t)\leq c \{a_1(t) + \sum^8_{i=4}a_i(t)\}.
\end{eqnarray}
By (3.14), (3.20) and (3.24)
\begin{eqnarray}
\Phi_\lambda (f(m_t)  )&\leq &1+ c\int^{t\wedge
\tau_2}_0\Phi_\lambda' (f(m_s)  )da_1(s)\nonumber\\
& +& c \sum^8_{i=4}\int^{t\wedge \tau_2}_0\Phi_\lambda'
(f(m_s)  ) da_i(s)   \nonumber\\
& +& 2\int^{t\wedge \tau_2}_0\Phi_\lambda'' (f(m_s) ){\bf trace}\{
m_s\cdot \widetilde{\sigma_s}\cdot \widetilde{\sigma_s}^T
m_s^T\}ds.\nonumber\\
\end{eqnarray}
So, taking mathematical expectation at both sides of (3.25), we have
\begin{eqnarray}
{\bf E}\{\Phi_\lambda (f(m_t)  )\}&\leq &1+ c\sum^8_{i=4} {\bf
E}\big \{ \int^{t\wedge \tau_2}_0\Phi_\lambda'
(f(m_s)  ) da_i(s) \big \}  \nonumber\\
& +& 2{\bf E}\big \{\int^{t\wedge \tau_2}_0\Phi_\lambda'' (f(m_s)
){\bf trace}\{ m_s\cdot \widetilde{\sigma_s}\cdot
\widetilde{\sigma_s}^T
m_s^T\}ds\big \}\nonumber\\
&:=& 1+ J_1(t) + J_2(t) + J_3(t) + J_4(t) + J_5(t)
+J_6(t).\nonumber\\
\end{eqnarray}
Now we estimate the terms $ J_i(t)$, $i=1, \cdots, 6$. \\

\noindent By (3.11), (3.12) and $N_t$ is bounded we get that
\begin{eqnarray}
J_1(t)&=& c{\bf E}\big \{ \int^{t\wedge \tau_2}_0\Phi_\lambda'
(f(m_s))N_s\big [2m_s\cdot \widetilde{b}_s + {\bf trace
}\{\widetilde{\sigma}_s\widetilde{\sigma}_s^T\} \big ]ds\big
\}\nonumber\\
& \leq & c^2 {\bf E}\big \{ \int^{t\wedge \tau_2}_0\frac{\lambda
\Phi_\lambda (f(m_s))}{\delta_1^2 + |m_s|^2}\big [
(3+\varepsilon^2)(\delta_1^2
+|m_s|^2 )\big ]ds\big \}\nonumber\\
& \leq & \lambda (3+\varepsilon )c^2{\bf E} \big \{\int^{t\wedge
\tau_2 }_0\Phi_\lambda (f(m_s))   \big \}ds.
\end{eqnarray}
\begin{eqnarray}
|J_2(t)|&=&\big | -\frac{2c}{\alpha} {\bf E}\big \{ \int^{t\wedge
\tau_2}_0\Phi_\lambda' (f(m_s))f(m_s)N_s \big [ (\triangledown
\phi^T \cdot {b} )(X^\varepsilon_{s\wedge \tau_2}(x)
)\nonumber\\
&&\qquad   \qquad \quad + (\triangledown \phi^T \cdot {b} )(
X^\varepsilon_n(s\wedge\tau_2,x) )         \big ]\big\} ds\big
|\nonumber\\
&\leq &  \lambda c  {\bf E}\big \{ \int^{t\wedge \tau_2}_0\frac{
\Phi_\lambda
(f(m_s))}{\delta_1^2 + |m_s|^2}|m_s|^2 ds\big \}\nonumber\\
& \leq & \lambda c  {\bf E}\big \{ \int^{t\wedge \tau_2}_0
\Phi_\lambda (f(m_s)) ds\big \}
\end{eqnarray}
because $\triangledown \phi  $, $ b$ and $N$ are bounded.\\

\noindent Similarly,
\begin{eqnarray}
|J_3(t)|&=&\big |\frac{\varepsilon}{\alpha}  {\bf E}\big \{
\int^{t\wedge \tau_2}_0\Phi_\lambda' (f(m_s)) N_sf(m_s) {\bf
trace}\big \{(\triangledown^2\phi \cdot \sigma \cdot \sigma^T )(
X^\varepsilon_{s\wedge \tau_2}(x))\nonumber\\
&& \qquad \qquad \quad +(\triangledown^2 \phi\cdot \sigma \cdot
\sigma^T )(X^\varepsilon_n(s\wedge\tau_2,x)
    ) \big\}ds\}\big |\nonumber\\
&&\leq  \lambda c \varepsilon {\bf E}\big \{ \int^{t\wedge \tau_2}_0
\Phi_\lambda (f(m_s)) ds\big \} .
\end{eqnarray}
\begin{eqnarray}
|J_4(t)|&=&\frac{2\varepsilon}{\alpha^2}\big | {\bf E}\big \{
\int^{t\wedge \tau_2}_0\Phi_\lambda' (f(m_s))N_s {\bf trace}\{\big (
\triangledown \phi^T \cdot {\sigma} \cdot \triangledown \phi \cdot
{\sigma^T} \big
)( X^\varepsilon_{s\wedge \tau_2}(x)    )\nonumber\\
&&\qquad \qquad \quad + \big ( \triangledown \phi^T \cdot {\sigma}
\cdot \triangledown \phi \cdot {\sigma^T} \big )(
X^\varepsilon_n(s\wedge\tau_2,x)    )
\big \}ds\big |\nonumber\\
&&\leq  \lambda c \varepsilon {\bf E}\big \{ \int^{t\wedge \tau_2}_0
\Phi_\lambda (f(m_s)) ds\big \} .
\end{eqnarray}
Using (3.12), and $ \sigma $ and $ \triangledown \phi  $ are bounded
we have
\begin{eqnarray}
|J_5(t)|&=& \frac{4\sqrt{\varepsilon}}{\alpha} \big | {\bf E}\big \{
\int^{t\wedge \tau_2}_0\Phi_\lambda' (f(m_s)) N_s{\bf trace}\big
\{\widetilde{\sigma_s}^Tm_s^T \big [ (\triangledown \phi^T \cdot
{\sigma} )(X^\varepsilon_{s\wedge \tau_2}(x)
)\nonumber\\
&&\qquad \qquad \quad +( \triangledown \phi^T \cdot {\sigma} )(
X^\varepsilon_n(s\wedge\tau_2,x) ) \big ]\big \}ds\big |\nonumber\\
&&\leq \lambda c{\varepsilon}{\bf E}\big \{\int^{t\wedge
\tau_2}_0\frac{ \Phi_\lambda
(f(m_s))}{\delta_1^2 + |m_s|^2} \sqrt{\delta_1^2 + |m_s|^2 }|m_s|     ds\big \}\nonumber\\
&& \leq \lambda c  \varepsilon {\bf E}\big \{ \int^{t\wedge
\tau_2}_0 \Phi_\lambda (f(m_s)) ds\big \}.
\end{eqnarray}
Since $ \Phi_\lambda'' (x) = \frac{\lambda  ( \lambda
-1)\Phi_\lambda (x) }{(
 \delta_1^2  + x )^2} $, $\forall  x\geq 0   $, by (3.12), we get
 that
\begin{eqnarray}
|J_6(t)|&=&\big | 2{\bf E}\big \{\int^{t\wedge
\tau_2}_0\Phi_\lambda'' (f(m_s) ){\bf trace}\{ m_s\cdot
\widetilde{\sigma_s}\cdot \widetilde{\sigma_s}^T m_s^T\}ds\big
\}\big |\nonumber\\
&& \leq \lambda ( \lambda -1)  c{\varepsilon}{\bf E}\big
\{\int^{t\wedge \tau_2}_0\frac{ \Phi_\lambda
(f(m_s))}{(\delta_1^2 + |m_s|^2)^2} ({\delta_1^2 + |m_s|^2 })|m_s|^2     ds\big \}\nonumber\\
&&\leq  \lambda  ( \lambda -1)c    \varepsilon {\bf E}\big \{
\int^{t\wedge \tau_2}_0 \Phi_\lambda (f(m_s)) ds\big \}.
\end{eqnarray}
Therefore, putting ( 3.27)-(3.32) and (3.26) together, we have
\begin{eqnarray*}
{\bf E} \{ \Phi_\lambda (f(m_t))\} \leq 1 + c[\lambda + \lambda
\varepsilon  + \lambda^2 \varepsilon ]\int^t_0 {\bf E} \{
\Phi_\lambda (f(m_{s\wedge \tau_2 }))\}ds,
\end{eqnarray*}
which, by Gronwall lemma, implies that
\begin{eqnarray}
{\bf E} \{ \Phi_\lambda (f(m_t))\} \leq  \exp\big \{ c[\lambda +
\lambda \varepsilon  + \lambda^2 \varepsilon ]t   \big \}.
\end{eqnarray}
Letting $t=1$ in the last inequality  we have
$$ {\bf P} \{ \tau_1\geq 1, \tau\leq 1  \} \Phi_\lambda (\delta^2    )
 \leq   \exp\big \{ c[\lambda +
\lambda \varepsilon  + \lambda^2 \varepsilon ] \big \}.
$$
So
\begin{eqnarray}
 \varepsilon \log {\bf P} \{ \tau_1\geq 1, \tau\leq 1  \} \leq
c[\lambda \varepsilon + \lambda \varepsilon^2  + \lambda^2
\varepsilon^2 ] -\lambda\varepsilon
\int^{\delta^2}_0\frac{ds}{\delta_1^2 +s}.
\end{eqnarray}
Taking $ \lambda\varepsilon =1 $, then making $\varepsilon $ tend to
$ +\infty $, finally letting $ \delta_1 $ tend to $0$ in (3.34), we
have
\begin{eqnarray}
\limsup_{\varepsilon \rightarrow 0} \varepsilon \log {\bf P} \{
\tau_1\geq 1, \tau\leq 1 \}=-\infty.
\end{eqnarray}
Since
\begin{eqnarray*}
{\bf P} \big\{ \sup_{t\in [0,1]}|
X^\varepsilon_t(x)-X^\varepsilon_n(t,x)|\geq \delta \big \}\leq {\bf
P}\{\tau_1\leq 1\} + {\bf P}\{\tau_1\geq 1, \tau_2\leq 1\},
\end{eqnarray*}
we immediately deduce  from (3.7) and (3.35)  that
\begin{eqnarray}
\lim_{n\rightarrow \infty }\limsup_{ \varepsilon\rightarrow 0  }
\varepsilon \log {\bf P} \big\{ \sup_{t\in [0,1]}|
X^\varepsilon_t(x)-X^\varepsilon_n(t,x)|\geq \delta \big \}=-\infty
\end{eqnarray}
for any $ \delta>0  . $ Therefore, we complete the proof of Theorem
3.1. \ \qquad $ \Box $
 \vskip 0.3cm
\setcounter{equation}{0}
\section{ \bf{\small  Large deviations principles on   SDE(1.4)with reflection}  }
\noindent
The main purpose of this section is to Theorem 1.1.\\

\noindent {\bf Proof of Theorem 1.1.}\  By (2.22) and (3.1), for
$n\geq 1 $ we have
$$ F_n(\sqrt{\varepsilon}B)=X^\varepsilon_n.$$
Using the schilder theorem (see Theorem 1.3.27 in \cite{s11} ), $
\{\sqrt{\varepsilon}B,\varepsilon>0\}$
 satisfies the large deviations
principles on $ E_1 $ with good rate function $ I_1(\cdot ) $
defined by (1.6). Because of Theorem 2.3 $F_n( \cdot ) $ is a
continuous map from $ E_1$ to $E$, therefore, by the contraction
principle(see Theorem 4.2.1 \cite{s10}), $\{ X^\varepsilon_n,
\varepsilon
>0 \}$ satisfies the large deviations principles on $E$ with good
rate function $I^x_{2,n}$   defined by
$$ I^x_{2,n}(f)=\inf\big \{I_1( \psi ): f=z_n^\psi\ \mbox{and } \
\psi\in L^2 ([0,1], \Re^d)    \big\}        $$ for $f\in E$, where
$z^\psi_n $ is solution of Eq.(2.14). By Theorem 2.2 and Theorem 3.1
as well as the generalized contraction principle for large
deviations principles (see Theorem 4.2.23 in \cite{s10}), the family
$\{ X^\varepsilon, \varepsilon >0  \}$ satisfies the large
deviations principles on $ E$ with the rate function $ I^x_2(\cdot )
$ defined by (1.8). So we complete the proof. \qquad $\Box$\vskip
19pt

Using the same way as in Theorem 1.1, we can prove the following and
we omit its proof here.
\begin{Them}
Assume the conditions of Theorem 1.1. $\{X^\varepsilon(t, y)  \}$ is
solution of SDE(1.4) corresponding to the initial condition $X_0=y$.
Then for any
 open set $G$ and any closed set $F$ of $E$ we have
\begin{eqnarray}
&& \liminf_{\stackrel{\varepsilon \rightarrow 0}{ y\rightarrow
x}}\varepsilon \log {\bf P}\{
X^\varepsilon(y)\in G\}\geq -\inf_{f\in G}\{I^x_2(f)\},\\
&& \limsup_{\stackrel{\varepsilon \rightarrow 0}{ y\rightarrow x}
}\varepsilon \log {\bf P}\{ X^\varepsilon(y)\in F\}\leq -\inf_{f\in
F}\{I^x_2(f)\}.
\end{eqnarray}
\end{Them}
As a direct consequence of Theorem 4.1 we also have the following.
\begin{Prop}
Assume the conditions of Theorem 1.1. Then  for any compact set
$K\subset \bar{ \mathcal{O}  }$, any open set $G$ and any closed set
$F$ of $E$ we have
\begin{eqnarray}
&& \liminf_{\varepsilon \rightarrow 0}\varepsilon \log \inf_{y\in K
} {\bf P}\{
X^\varepsilon(y)\in G\}\geq -\sup_{y\in K}\inf_{{f\in G}{}}\{I^y_2(f)\},\\
&& \limsup_{\varepsilon \rightarrow 0
 }\varepsilon \log \sup_{y\in
K}{\bf P}\{ X^\varepsilon(y)\in F\}\leq -\inf_{\stackrel{f\in
G}{y\in K}}\{I^y_2(f)\}.
\end{eqnarray}
\end{Prop}
\vskip 0.3cm \setcounter{equation}{0}
\section{{\bf{\small  Large deviations principles on SDE(1.4)with
reflection and anticipated initial conditions} }} \noindent Let $
X^\varepsilon_0 $ be any random variable, which may not be adapted
to $ \mathcal{F}_0 $, $ X^\varepsilon_t(x)$ be a solution of
SDE.(1.4) corresponding to initial data $x$.  The author proved in
\cite{s6} that the composition  $ X^\varepsilon_t( X^\varepsilon_0
)$ is also a solution of SDE ( 1.4) corresponding to anticipated
initial conditions  $ X^\varepsilon_0 $. Therefore the integral in
this equation should be  anticipated in Malliavin sense(see
\cite{s12}). The main purpose of this section is to prove Theorem
1.2 via Theorem 1.1 and Proposition 4.1. The proof is the following.
\vskip 0.3cm \noindent {\bf Proof of Theorem 1.2.}\ {\sl We first
prove that the large deviations lower bound holds for $
\{X^\varepsilon_t( X^\varepsilon_0 ), \varepsilon >0   \}$ with rate
function $ I^{x_0}(\cdot )$ on $E=C([0,1], \bar {\mathcal{O} })$ if
the conditions of Theorem 1.2 are satisfied.}

 Let $ G\subset
E $ be an open set. Assume that $g\in G $ with $ I^{x_0}(g) <
+\infty $. By the definition of $I^{x_0}(g) $, there exist $
\delta_2 >0$, $ \psi \in L^2([0,1]; \Re^d ) $ and $ z^\psi(y)$, a
solution of Eq.(1.5) corresponding to initial $y$, such that $ \
|y-x_0| < \delta_2 $ and
\begin{eqnarray}
I_2^y(g)=I_2^y(z^\psi(y))=I_1( \psi )\leq I^{x_0}(g) < + \infty .
\end{eqnarray}
Using the same way as in (2.8), we can prove that for any $ y_1, y_2
\in \Re^d $
\begin{eqnarray}
\sup_{t\in [0,1]}\big \{| z^\psi(y_1)-z^\psi(y_2)    |^2\big \} \leq
c|y_1-y_2|^2\exp\{c(1+ \|\psi\|_{L_2})\}
\end{eqnarray}
for some constants $c=c(\alpha, C_0, O)$. Taking $ \delta_3 >0 $
such that
\begin{eqnarray}
B_{\delta_3 }(y)= \big \{f: \| f- z^\psi(y) \|_{E}\leq \delta_3 \big
\},
\end{eqnarray}
we have
\begin{eqnarray}
{\bf P}\big \{ X^\varepsilon_\cdot( X^\varepsilon_0 )\in G \big \}
&\geq & {\bf P}\big \{ \|X^\varepsilon_\cdot( X^\varepsilon_0
)-z^\psi(y)\|_{E}\leq \delta_3, \ | X^\varepsilon_0 -y|\leq
2\delta_3  \big \}\nonumber\\
&\geq & {\bf P}\big \{ \sup_{|x-y|\leq
2\delta_3}\|X^\varepsilon_\cdot( x
 )-z^\psi(y)\|_{E}\leq \delta_3, \ |
X^\varepsilon_0 -y|\leq 2\delta_3  \big \}\nonumber\\
&&\mbox{(  by (5.2))  }\nonumber\\
&\geq & {\bf P}\big \{ \sup_{|x-y|\leq
2\delta_3}\|X^\varepsilon_\cdot( x
 )-z^\psi(x)\|_{E}\leq \frac{\delta_3 }{2}, \ |
X^\varepsilon_0 -y|\leq 2\delta_3 \big \}\nonumber\\
&= & {\bf P}\big \{ \sup_{|x-y|\leq 2\delta_3}\|X^\varepsilon_\cdot(
x )-z^\psi(x)\|_{E}\leq \frac{\delta_3 }{2}\big \}\nonumber\\
&&\qquad  \qquad \qquad -{\bf P}\big \{
| X^\varepsilon_0 -y| > 2\delta_3 \big \}\nonumber\\
&\geq & \inf_{|x-y|\leq 2 \delta_3}{\bf P}\big \{
\|X^\varepsilon_\cdot( x )-z^\psi(x)\|_{E} < \frac{\delta_3
}{2}\big \}\nonumber\\
&& \qquad \qquad -{\bf P}\big \{ | X^\varepsilon_0 -y| > 2\delta_3
\big \}.
\end{eqnarray}
Using (4.3), (5.1) and Theorem 1.1 we have
\begin{eqnarray}
&&\liminf_{\varepsilon \rightarrow 0}\varepsilon \log
\inf_{|x-y|\leq 2 \delta_3}{\bf P}\big \{ \|X^\varepsilon_\cdot( x
)-z^\psi(x)\|_{E} < \frac{\delta_3 }{2}\big \}\nonumber\\
&&\geq \inf_{ |x-y| \leq 2\delta_3}\big [ -\inf_{f\in
B_{\frac{\delta_3}{2}} ( z^\psi(y)  )} \big \{I_2^y(f)     \big
\}\big ]\nonumber\\
&&\geq \inf_{ |x-y| \leq 2\delta_3}\big [ - I_2^y( z^\psi(y) ) \big ]\nonumber\\
&&\geq \inf_{ |x-y| \leq 2\delta_3}\big [ - I_1( \psi) ) \big ]=-
I_1( \psi) \geq - I^{x_0}(g).
\end{eqnarray}
Therefore,  putting (1.14), (5.4) and (5.5) together, we have
\begin{eqnarray*}
\liminf_{\varepsilon \rightarrow 0}\varepsilon \log {\bf P}\big \{
X^\varepsilon_\cdot( X^\varepsilon_0 )\in G \big \}\geq -
I^{x_0}(g).
\end{eqnarray*}
Since $ g\in G $ is arbitrary, we have
\begin{eqnarray*}
\liminf_{\varepsilon \rightarrow 0}\varepsilon \log {\bf P}\big \{
X^\varepsilon_\cdot( X^\varepsilon_0 )\in G \big \}\geq - \inf_{
f\in G} \big \{  I^{x_0}(g)\big \}.
\end{eqnarray*}
So we complete the proof of the large deviations lower
bound.\vskip19pt

{\sl  Next we prove that  the large deviations upper bound holds for
$ \{X^\varepsilon_t( X^\varepsilon_0 ), \varepsilon >0   \}$ with
rate function $ I^{x_0}(\cdot )$ on $E$ if the conditions of Theorem
1.2 are satisfied.}

 Let $ F \subset E  $ be a closed set. If $ \inf_{f\in F } \big \{I^{x_0}(f)   \} =0
$, then (1.16) is trivial. We need only to prove that for any $a>0 $
with $ a < \inf_{f\in F } \big \{I^{x_0}(f)   \}$
\begin{eqnarray}
\limsup_{\varepsilon \rightarrow 0}\varepsilon \log {\bf P}\big \{
 X^\varepsilon_\cdot( X^\varepsilon_0 )\in F\big \} \leq -a.
\end{eqnarray}
Now we assume that $ \exists \ a>0  $ such that
$$   \inf_{f\in F } \big \{I^{x_0}(f)   \}> a.  $$
Then there exists $ \delta_4 >0 $ and $ y_0 \in B(x_0, \delta_4 )=\{
x: |x-x_0|\leq \delta_4 \}$ such that
\begin{eqnarray}
\inf_{f\in F}\big\{ I^{y_0}_2(f)    \big \}>a .
\end{eqnarray}
Then the level set $ K_a(y_0):=\big \{ f: I^{y_0}_2(f)\leq a \big
\}=\big \{ z^\psi(y_0) : I^{y_0}_2(z^\psi(y_0))\leq a, \  \psi \in
L^2([0,1]; \Re^d     ) \big \},$ $z^\psi(y_0) $ is a solution of
Eq.(1.5) corresponding to initial data  $y_0$. Since $K_a(y_0)\cap
F=\emptyset $, $ \forall \  f=z^\psi(y_0)  \in   K_a(y_0)  $, $
\exists \ \delta_f >0 $ such that
$$ U_{\delta_f}(f)=\big \{ g: \|g-f\|_{E_2} < \delta_f    \} \cap F=\emptyset. $$
Noting that $\bigcup\limits_{f\in  K_a(y_0) } U_{\delta_f}(f)\supset
K_a(y_0) $ and $ K_a(y_0) $ is compact, $\exists    $ $m>0 $ and
$f_j \in  K_a(y_0) $ such that
\begin{eqnarray}
{\bf U}:=\bigcup_{j=1}^{m} U_{\delta_{f_j}}(f_j)\supset K_a(y_0).
\end{eqnarray}
So
\begin{eqnarray}
{\bf U}\cap F =\emptyset.
\end{eqnarray}
Hence
\begin{eqnarray}
&&{\bf P}\big \{ X^\varepsilon_\cdot( X^\varepsilon_0 )\in F \big
\}\nonumber \\
&&\leq {\bf P}\big\{ X^\varepsilon_\cdot( X^\varepsilon_0 )\in F,
|X^\varepsilon -y_0|\leq 2 \delta_4 \big\}+{\bf P}\big
\{|X^\varepsilon_0-x_0|\geq \delta_4  \}\nonumber\\
&&:= \Box_1 + \Box_2.
\end{eqnarray}
By (5.9)
\begin{eqnarray}
\Box_1 &\leq  &\sup_{|y-y_0|\leq 2\delta_4}{\bf P}\big \{
X^\varepsilon_\cdot( y )\in F \big \}\nonumber\\
&\leq &\sup_{|y-y_0|\leq 2\delta_4}{\bf P}\big \{
X^\varepsilon_\cdot( y )\in {\bf U}^c \big \}\nonumber\\
&\leq & \sup_{|y-y_0|\leq 2\delta_4}{\bf P}\big \{
X^\varepsilon_\cdot( y_0 )\in {\bf U}^c \big \}
\nonumber\\&&+\sup_{|y-y_0|\leq 2\delta_4}\limsup_{k\rightarrow
\infty }{\bf P}\big \{|X^\varepsilon_\cdot( y
)-X^\varepsilon_\cdot( y_0 )|\geq \frac{1}{k} \big \}\nonumber\\
&\leq & \sup_{|y-y_0|\leq 2\delta_4}{\bf P}\big \{
X^\varepsilon_\cdot( y_0 )\in {\bf U}^c \big \}
\nonumber\\
&&+\limsup_{k\rightarrow \infty }\sup_{|y-y_0|\leq 2\delta_4}{\bf
P}\big \{|X^\varepsilon_\cdot( y
)-X^\varepsilon_\cdot( y_0 )|\geq \frac{1}{k} \big \}\nonumber\\
 &:=&\Box_{11} + \limsup_{k\rightarrow \infty }\Box_{12}.
\end{eqnarray}
Using Theorem 1.1 and ${\bf U}^c  $ is a closed set, we have
\begin{eqnarray}
\limsup_{\varepsilon\rightarrow 0 }\varepsilon \log \Box_{11} \leq -
\inf_{f \in {\bf U}^c  }\{ I^{y_0}_2(f) \} \leq -a
\end{eqnarray}
because ${\bf U}^c \subset K_a(y_0)^c $. \\
We claim that
\begin{eqnarray}
\limsup_{k\rightarrow \infty }\limsup_{\delta_4\rightarrow 0}
\limsup_{\varepsilon\rightarrow 0} \varepsilon \log \Box_{12} =
-\infty.
\end{eqnarray}
To prove (5.13) we define for $k\geq 1$
\begin{eqnarray}
&&\tau =\inf\{t: |X^\varepsilon_t(y)- X^\varepsilon_t(y_0)  |\geq
\frac{1}{k}   \}\wedge 1,\\
&&m_t = X^\varepsilon_{t\wedge \tau}(y)- X^\varepsilon_{t\wedge
\tau}(y_0),\\
&&b_1(t) =b(X^\varepsilon_{t\wedge \tau}(y))-b(
X^\varepsilon_{t\wedge
\tau}(y_0)),\\
&&\sigma_1(t) =\sqrt{\varepsilon}\big [\sigma(X^\varepsilon_{t\wedge
\tau}(y))-\sigma( X^\varepsilon_{t\wedge \tau}(y_0))\big ].
\end{eqnarray}
Then
\begin{eqnarray}
m_t&=&y-y_0 + \int^t_0 b_1(s)ds + \int^t_0 b_1(s)dB(s)\nonumber\\
&& -[L^\varepsilon_{t\wedge \tau}(y)- L^\varepsilon_{t\wedge
\tau}(y_0)].\nonumber\\
\end{eqnarray}
Define functions $ \Psi $  and $\Phi_\lambda  $ on $[0, +\infty )$
by
$$ \Psi (x)= \int^x_0 \frac{ds}{\rho^2+ s} \ \mbox{ and } \
\Phi_\lambda(x)=\exp\big \{ \lambda \Psi(x)  \big \} \ \mbox{for
$\lambda >0$  and $ \rho >0 $.} $$ By the same way as in (3.36), we
obtain that
\begin{eqnarray}
\limsup_{\delta_4\rightarrow 0} \limsup_{\varepsilon\rightarrow 0}
\varepsilon \log \Box_{12} \leq -
\int^{\frac{1}{k^2}}_0\frac{ds}{\rho^2 +s} + c
\end{eqnarray}
for some constants $c$ and any $k\geq 1$. Letting $\rho \rightarrow
0$ we have for any $k\geq 1$
\begin{eqnarray}
\limsup_{\delta_4\rightarrow 0} \limsup_{\varepsilon\rightarrow 0}
\varepsilon \log \Box_{12} = - \infty.
\end{eqnarray}
This implies claim (5.13).\\
Therefore, we deduce from that (5.11) to (5.13) that
\begin{eqnarray}
\limsup_{\varepsilon \rightarrow 0} \varepsilon \log \Box_1 &  \leq
& \big \{\limsup_{\varepsilon\rightarrow 0 }\varepsilon \log
\Box_{11}\big\} \nonumber \\
&&\vee \big \{ \limsup_{k\rightarrow \infty
}\limsup_{\delta_4\rightarrow 0} \limsup_{\varepsilon\rightarrow 0}
\varepsilon \log \Box_{12}\big \}\nonumber\\
&\leq &-a.
\end{eqnarray}
So we know from (1.13),  (5.10) and (5.21) that
\begin{eqnarray}
\limsup_{\varepsilon \rightarrow 0} \varepsilon \log {\bf P}\big \{
X^\varepsilon_\cdot( X^\varepsilon_0 )\in F \big \}& \leq &
\{\limsup_{\varepsilon \rightarrow 0} \varepsilon \log \Box_1\}
\nonumber \\&&  \vee \{\limsup_{\varepsilon\rightarrow
0}\varepsilon\log {\bf P}\{|X^\varepsilon_0-x_0|>\delta_4  \}
\}\nonumber\\ &\leq & -a,
\end{eqnarray}
that is, (5.6) holds. Thus we complete the proof of of the large
deviations upper bound. \vskip19pt

 {\sl  Finally, we prove that $ I^x(\cdot )$
is a good rate function on $E$.}

 Since for any $ a \in [0, \infty )
$
$$\{f: I^x(f)>a \}= \bigcup_{n=1} \bigcup_{y: |y-x_0|\leq \frac{1}{n}}
  \{f: I_2^y(f)>a  \} $$
and $I^y_2(\cdot ) $ is a good rate function, we know that $ \{f:
I^x(f)>a \}    $ is an open set and $ \{f: I_2^y(f) \leq a\}      $
is compact set. So $I^x $ is a lower semicontinuous function.
Because $ \{f: I^x(f)\leq a \}\subset \{f: I_2^y(f) \leq a\}$ for
$y\in B(x_0, 1)$, $ \{f: I^x(f)\leq a \}$ is also compact, so  $
I^x(\cdot )$ is also a good rate function on $E$. Thus we complete
the
proof of Theorem 1.2. \qquad $\Box $.\\

{\bf Acknowledgements.} This work is supported by NSFC and SRF for
ROCS, SEM. The author would like to thank both for their generous
financial support.

\end{document}